\documentclass[14pt]{book}

\begin {document}
\newtheorem{theorem}{Theorem}[section]

\newtheorem{corollary}[theorem]{Corollary}
\newtheorem{lemma}[theorem]{Lemma}
\newtheorem{fact}{Example}

\newenvironment{proof}{
\par
\noindent {\bf Proof.}\rm}

\title{Simple decompositions of simple Lie superalgebras}
\author{T.V.Tvalavadze
\thanks{\it supported by NSERC Grant 227060-00}\\
 \small Department of Mathematics and Statistics\\
\small Memorial University of Newfoundland\\
\small St. John's, NL, CANADA.}
\date{}
\maketitle

\sloppy

\chapter{Introduction}

In this paper we consider Lie superalgebras decomposable as the
sum of two proper  subalgebras. Any of these algebras has the form
of the vector space sum $L=A+B$ where $A$ and $B$ are proper
simple subalgebras which need not be ideals of  $L$, and the sum
need not be direct.

The main result  of this paper is the following
\begin{theorem}

Let  $ S = {osp}(m,2n)$ be a  Lie superalgebra such that $S=K+L$
where $K$, $L$ are two proper basic simple  subalgebras.  Then $m$
is even, $m=2k$ and $K \cong osp(2k-1,2n)$, $L \cong sl(k,n)$.
\end{theorem}

In this paper all suablebras are $Z_2$-graded.

 {\section{Decompositions of $osp(m,2n)$ as the sum of basic Lie
subalgebras }

\subsection{ Preliminaries  }

We use the following technical Lemmas:

\begin{lemma}

Let $L \cong sl(m,n)$, $L=L_0 \oplus L_1$ where $L_0$ is an even
part of $L$, $L_1$ is an odd part of $L$. Then the following
properties hold:

\noindent (a)  $L_0 = I_1 \oplus I_2 \oplus U $, where $I_1\cong
sl(m)$, $I_2\cong sl(n)$ and  $U$ is either one dimensional Lie
algebra if $m=n$ or zero element.

\noindent (b) $L_0$-module $L_1$ is a direct sum of two
irreducible $L_0$-modules of the dimension~$mn$ with highest
weights $(\lambda,\mu)$ and $(\mu,\lambda)$, where $\lambda =(1,0,
\ldots, 0 )$ and  $\mu =(0, \ldots, 0,1 )$.

\noindent (c) $[L_1,L_1]=L_0$

\noindent (d) $[I_1,L_1]=L_1$ and $[I_2,L_1]=L_1$

\noindent (e) $I_1$-module $L_1$ is a direct sum of $2n$
irreducible $I_1$-modules of  dimension $m$ and $I_2$-module $L_1$
is a direct sum of $2m$ irreducible $I_2$-modules of
dimension~$n$.

\end{lemma}

\begin{lemma}

Let $L \cong osp(m,2n)$. Then

\noindent (a)  $L_0 = I_1 \oplus I_2$, where $I_1\cong o(m)$,
$I_2\cong sp(2n)$ .

\noindent (b) $L_0$-module $L_1$ is a irreducible $L_0$-modules of
 dimension $2mn$.

\noindent (c) $[L_1,L_1]=L_0$

\noindent (d) $[I_1,L_1]=L_1$ and $[I_2,L_1]=L_1$

\noindent (e) $I_1$-module $L_1$ is a direct sum of $2n$
irreducible $I_1$-modules of  dimension $m$ and $I_2$-module $L_1$
is a direct sum of $m$ irreducible $I_2$-modules of dimension
$2n$.
\end{lemma}

The proof  of these Lemmas can be found in \cite{Kac}.

The following two Lemmas give the decomposition of simple Lie
algebra as the sum of simple subalgebras. They were found by
 Onishchik (see  \cite{ON}). These matrix forms can be found in
 \cite{BTT}.
\begin{lemma}
\par\medskip
 Let $o(2n)$ be decomposable into the sum of two subalgebras isomorphic
 to $o(2n-1)$ and $sl(n)$. Then there exists a basis of $F^{2n}$ such that this decomposition
takes the following matrix form:
$$S=N+M, \eqno(1) $$
where  $S \cong o(2n)$ consists of the matrices:
$$
 \qquad
    \left( \begin{array}{c|c}
   A_{11} &  A_{12}     \\
  \hline
   A_{21} &  A_{22}    \\
  \end{array}\right)
\eqno(2)
$$ where $A_{12}, A_{21}$ are skew-symmetric  matrices of order $n$ and
$A_{11}, A_{22}$ are matrices of order $n$ such that
$A_{22}=-A_{11}^{t}.$

\noindent  The first subalgebra $N\cong o(2n-1)$ consists of the
matrices:
$$
    \left( \begin{array}{c|ccc|c|ccc}
     0&  y_1& \ldots&   y_{n-1}&      0& x_1&   \ldots & x_{n-1}  \\
     \hline
   x_1&    &       &      &   -x_1&    &           &     \\
\vdots&    &  A'_{11}  &      & \vdots&    &   A'_{12}     &     \\
   x_{n-1}&    &       &      &   -x_{n-1}&    &           &    \\
     \hline
     0& -y_1& \ldots&  -y_{n-1}&      0&-x_1&   \ldots & -x_{n-1}   \\
     \hline
    y_1&   &       &      &   -y_1&    &           &     \\
 \vdots&   &   A'_{21} &      & \vdots&    &  A'_{22}      &        \\
    y_{n-1}&   &       &      &   -y_{n-1}&    &           &      \\
     \end{array}\right)
           \eqno (3)
$$ where $A'_{12}, A'_{21}$ are skew-symmetric  matrices of order $n-1$ and
$A'_{11}, A'_{22}$ are matrices of order $n-1$ such that
$A'_{22}=-{A'}_{11}^{t}.$

\noindent The second subalgebra  $M\cong
 sl(n)$ consists of the matrices:
$$   \left( \begin{array}{c|c}
    A_{1} &  0     \\
  \hline
    0  &  A_{2}    \\
 \end{array}\right)
$$
\noindent where $A_1$, $A_2$ are matrices of order $n$ with zero
trace such that $A_2=-A_1^{t}.$

\end{lemma}

\begin{lemma}

 Let $o(4n)$ be decomposable into the sum of two subalgebras isomorphic
 to $o(4n-1)$ and $sp(2n)$. Then there exists a basis of $F^{4n}$ such that this decomposition
takes the following matrix form:
$$S=N+M \eqno(4) $$
where  $S \cong o(4n)$ consists of the matrices of the form (2)
where $A_{11}, A_{12}, A_{21}, A_{22}$ are of the order $2n$.

\noindent The first subalgebra $N\cong o(4n-1)$ has the  form (3),
where  $A_1, B_1, C_1, D_1$ are of the order $2n-1$.

\noindent  The second subalgebra $M\cong  sp(2n)$ consists of the
matrices:
$$
   \left( \begin{array}{c|c}
   Y &  0     \\
  \hline
    0  & -Y^t    \\
 \end{array}\right)
$$
where $Y$ is the set of matrices
$$
    \left( \begin{array}{c|c}
   A & B      \\
  \hline
   C &  D    \\
  \end{array}\right)
$$
where $B, C$ are skew-symmetric matrices of order $n$ and
$D=-A^{t}$.

\end{lemma}

\subsection{ Description of $L_0$-modules and $K_0$-modules }

 Let  $ S = {osp}(m,2n)$ be a Lie superalgebra such that $S=K+L$ where $K$,
$L$ are two proper basic simple subalgebras. We consider ${S}
\cong osp(m,2n)$ as the subalgebra of $gl(m,2n)$. Then ${L}
\subset { S}$ is also a subalgebra of $gl(m,2n)$ and ${L}_0
\subset gl(m) \oplus gl(2n)$. Hence we have two  natural
representations $\rho_1$ and $\rho_2$ of ${ L}_0$ in vector spaces
$V$ and $W$ where $V$ is a column vector space of dimension $m$,
and $W$ is a column vector space of dimension $2n$. We will also
consider $V$ and $W$ as ${L}_0$-module such that
$$xv=\rho_1(x)(v)$$ and $$xw=\rho_2(x)(w),$$ for any $x\in L_0$,
$v \in V$, $w \in W$.

$L_0$-modules $V$ and $W$ are completely reducible because $L_0$
is a reductive Lie algebra. Let $V=V_1 \oplus \ldots \oplus V_r$
and $W=W_1 \oplus \ldots \oplus W_d,$ where $V_i$,$W_j$ are simple
$L_0$-modules.

Next we consider $L_0$-module $V=V_1 \oplus \ldots \oplus V_r$.
Let $I_1$ and $I_2$ be ideals of $L_0$. Next we define
 the following types of $L_0$-module $V_i$.

\noindent {\bf 1.} $L_0$-module $V_i$ is of {\it{type 1}} if
$L_0$-module $V_i$ is trivial.

\noindent {\bf 2.} $L_0$-module $V_i$ is of {\it{type 2}} if $I_2$
acts trivially on $V_i$ but $I_1(V_i)\neq \{0\}$.
 In this case, we consider $V_i$ as $I_1$-module.

\noindent {\bf 3.} $L_0$-module $V_i$ is of {\it{type 3}} if $I_1$
acts trivially on $V_i$ but $I_2(V_i)\neq \{0\}$. In this case, we
consider $V_i$ as $I_2$-module.

\noindent {\bf 4.} $L_0$-module $V_i$ is of {\it{type 4}} if
$I_1(V_i)\neq \{0\}$ and $I_2(V_i)\neq \{0\}$

$\quad$

Similarly $L_0$-module $W_j$ can also be  one of the following
types:

\noindent {\bf 1.} $L_0$-module $W_j$ is of {\it{type 1}} if
$L_0$-module $W_j$ is trivial.

\noindent {\bf 2.} $L_0$-module $W_j$ is of {\it{type 2}} if $I_2$
acts trivially on $W_j$ but $I_1(W_j)\neq \{0\}$.

\noindent {\bf 3.} $L_0$-module $W_j$ is of {\it{type 3}} if $I_1$
acts trivially on $W_j$ but $I_2(W_j)\neq \{0\}$.

\noindent {\bf 4.} $L_0$-module $W_j$ is of {\it{type 4}} if
$I_1(W_j)\neq \{0\}$ and $I_2(W_j)\neq \{0\}$

$\quad$

In a similar  manner we  define   types of $K_0$-modules.

Now we look at the decomposition $S=K+L$. We consider ${S} \cong
osp(m,2n)$ as a subalgebra of $gl(m,2n)$. Hence $S_0=K_0+L_0
\subset gl(m,2n)_0 $ and $S_1=K_1+L_1 \subset gl(m,2n)_1$. There
exists  an isomorphism between vector spaces $gl(m,2n)_1$ and $(V
\otimes W^{\ast}) \oplus (V^{\ast} \otimes W)$. We identify
$gl(m,2n)_1$ with $(V \otimes W^{\ast}) \oplus (V^{\ast} \otimes
W)$. Hence ${L}_0$-module $gl(m,2n)_1$ can be viewed as a direct
sum of two ${L}_0$-modules $V \otimes W^{\ast}$ and $V^{\ast}
\otimes W$ such that $$x(v \otimes f) = \rho_1(x)(v)\otimes f + v
\oplus \rho^{\ast}_2(x)(f)$$ and $$x(g \otimes w)
=\rho^{\ast}_1(x)(g)\otimes w  + g \oplus \rho_2(x)(w),$$ for any
$ x\in L_0$, $v \in V$, $w \in W$, $g \in V^{\ast}$, $f \in
W^{\ast}$ and $\rho^{\ast}_1$, $\rho^{\ast}_2$ are dual
representations for $\rho_1$, $\rho_2$.

Since $V=V_1 \oplus \ldots \oplus V_r$ and $W=W_1 \oplus \ldots
\oplus W_d$ where $V_i$,$W_j$ are simple $L_0$-modules, we can
express $L_0$-module $V^{\ast} \otimes W$ as a direct sum of
$L_0$-modules $V_i^{\ast} \otimes {W_j}$, $$V^{\ast} \otimes W =
\bigoplus_{i,j}(V_i^{\ast} \otimes {W_j}).$$

We denote a projection of  $V^{\ast} \otimes W$ onto $V_i^{\ast}
\otimes {W_j}$ as $\pi_{ij}$. If $\dim{ W_{j_0}}= \dim{ W_{j_1}}$
then there exists an isomorphism $\lambda_{j_0j_1}:{ W_{j_0}}
\mapsto { W_{j_1}}$. It extends to an isomorphism $\mu_{ij_0j_1}:{
V_i^{\ast} \otimes W_{j_0}} \mapsto { V_i^{\ast} \otimes
W_{j_1}}$. In this paper instead of the expression
$\pi_{ij_{0}}(x)= \lambda\mu_{ij_0j_1}(\pi_{ij_{1}}(x))$ for any
$x\in L_1$ we  write $\pi_{ij_{0}}(L_1)=
\lambda\pi_{ij_{1}}(L_1)$.

 We choose a basis in $V \oplus W$ from elements of subspaces
$V_i$, $i=1 \ldots r$ and $W_j$, $j=1 \ldots d$. Let us identify
the elements form $S$ and their matrix realizations in this basis.
Let $\dim{V_i}=m_i$ and $\dim{W_j}=n_i $. Then  we denote as
$\varepsilon'_i$  the natural isomorphism of $V_i$ onto the column
vector space $F^{m_i}$ and $\varepsilon''_i$  the natural
isomorphism of $W_i$ onto the column vector space $F^{n_j}$. We
extend $\varepsilon'_i$ and $\varepsilon''_j$ to an isomorphism
$\varepsilon_{ij}:{ V_i^{\ast} \otimes W_{j}} \mapsto { F^{m_i}
\otimes F^{n_j}}$. Let us define $\varrho_{ij} =
\varepsilon_{ij}\pi_{ij}: V^{\ast} \otimes W \mapsto F^{m_i}
\otimes F^{n_j} $

\subsection{ Properties of subalgebras ${osp}(p,2q)$ and $sl(s,l)$ in
the decomposition }

\begin{lemma}

 Let $ S = {osp}(m,2n)$ be a Lie superalgebra, and $S$ be decomposed into the sum of two proper simple
subalgebras ${ K}$ and ${ L}$ of the  type ${osp}(p,2q)$ and
${sl}(s,l)$, respectively. Then  $m$ is even, $m=2k$ for some $k$,
$p=2k-1$, $q=n$ and either $s=k$ or $l=k$.

\end{lemma}

\begin{proof}
By Lemma 2.2(a),  ${S}_0 = {o(m)} \oplus sp(2n)$. We define two
projections $\pi_1$ and $\pi_2$ of ${S}_0  $ into the ideals
$o(m)$ and $sp(2n)$, $\pi_1:{ S}_0\to{o(m)}$ and $\pi_2:{ S}_0\to
{sp(2n)}$. Since  ${ K}$ is isomorphic to ${osp}(p,2q)$,  ${ K}_0$
is isomorphic to ${o(p)} \oplus {sp(2q)}$. By Lemma 2.1(a), ${
L}_0$ is isomorphic to  ${sl(s)} \oplus {sl(l)} \oplus { U}$.
Since ${ K}_0$ and ${L}_0$ are reductive subalgebras, the
projections $\pi_1({ K}_0)$, $\pi_1({ L}_0)$, $\pi_2({ K}_0)$ and
$\pi_2({ L}_0)$ are also reductive as homomorphic images of
reductive algebras.

Since ${ S} = { K} + { L}$,  ${ S}_0$ is decomposable into the sum
of two subalgebras ${K}_0$ and ${ L}_0$, ${ S}_0 = { K}_0 + {
L}_0$. Therefore, $\pi_1({ S}_0) = \pi_1({ K}_0) + \pi_1({ L}_0)$
and $\pi_2({ S}_0) = \pi_2({ K}_0) + \pi_2({ L}_0)$, where
$\pi_1({\cal S}_0)={ o(m)}$ and $\pi_2({\cal S}_0)={ sp(2n)}$. We
have the decompositions of simple Lie algebras $o(m)$ and $sp(2n)$
into the sum of two reductive subalgebras.

By Onishchik's Theorem (see \cite{On}), ${sp(2n)}$ cannot be
decomposed into the sum of two proper reductive subalgebras. Hence
$sp(2n) = \pi_2({ K}_0) + \pi_2({ L}_0)$ is a trivial
decomposition and $\pi_2({ K}_0)= sp(2n)$. Since $K_0 = o(p)
\oplus sp(2q)$, it follows that $q=n$.

By  Onishchik's Theorem,  $o(m)$ has two decompositions into the
sum of two proper reductive subalgebras:

1. If $m=2k$ then  $o(2k)=o(2k-1)+ sl(k),$

2. If $m=4k$ then  $o(4k)=o(4k-1)+sp(2k).$

 The decomposition $o(m) = \pi_1({ K}_0) + \pi_1({ L}_0)$ cannot
 be of the second type, because $\pi_1({ L}_0)$ is not isomorphic to
 $sp(2k)$. Moreover, this decomposition cannot be trivial, because
  $\pi_1({ K}_0) = o(m)$ and $\pi_2({ K}_0) = sp(2n)$. Hence  $K_0$ coincides
 with $S_0$. This contradicts the fact that $K$ is a proper
 subalgebra of $S$.

 Therefore  $o(m) = \pi_1({ K}_0) + \pi_1({ L}_0)$
 is the decomposition of the first type and $m=2k$,  $\pi_1({ K}_0) \cong
o(2k-1)$, $\pi_1({ L}_0) \cong
 sl(k)$. Since $K_0 \cong o(p) \oplus sp(2q)$ and $L_0 \cong {sl(s)} \oplus {sl(l)} \oplus {U}$,
   it follows  that $p=2k-1$, $q=n$ and either $s=k$ or $l=k$.
\hfill $\fbox{}$

\end{proof}

$\quad$

Without any loss of generality, we assume that $L \cong sl(k,l)$.

\begin{corollary}

$L_0$-module $V$ is a direct sum of two $L_0$-modules of type 2,
$V=V_1 \oplus V_2$. Moreover $I_2$ acts trivially on $V$,
$I_1$-module $V_1$ is standard, and $I_1$-module $V_2$ is dual.
\end{corollary}

\begin{proof}
From  Lemma 2.5  in the decomposition $o(2k) = \pi_1({ K}_0) +
\pi_1({ L}_0)$ the first component is isomorphic to $o(2k-1)$ and
the second one is isomorphic to $sl(k)$. By Lemma 2.3, there exist
bases of $V$ such that the decomposition $o(2k)=o(2k-1)+ sl(k)$
takes the matrix form (1). Hence  $\pi_1({ L}_0)$ takes a form:
$$\left\{    \left(\begin{array}
                   { c  |  c    }
                     Y  & 0     \\ \hline
                     0  & -Y^{t}     \\
           \end{array}
   \right)\right\} , $$
   where  $Y$  is a  set of matrices of order \textit{k}  with a zero
   trace. Therefore we obtain   that $V=V_1 \oplus V_2$ and $I_2$ acts
trivially on $V$.  Moreover,  $I_1$-module $V_1$ is standard and
$I_1$-module $V_2$ is dual.

 \hfill $\fbox{}$

\end{proof}

$\quad$

\subsection{ Decompositions of  $osp(2k,2n)$ as
the sum of \\  $osp(2k-1,2n)$ and $sl(k,l)$ }

\begin{lemma}
Let $S=K+L $, where $S\cong osp(2k,2n)$, $K\cong osp(2k-1,2n)$,
$L\cong sl(k,l)$.

\noindent (a) There is no  $j_{0} \in \{1 \ldots d\}$ such that
for any $i \in \{1,2\}$, $\pi_{ij_{0}}(L_1)=0$.

\noindent (b) There are no $j_0,j_1 \in \{1 \ldots d\}$ such that
for some $\lambda \in F$, $\varrho_{ij_{0}}(x)=
\lambda\varrho_{ij_{1}}(x)$ where $x \in L_1$, $i \in \{1,2\}$.

\end{lemma}

\begin{proof}
By Corollary 2.6, $V=V_1 \oplus V_2$ as simple $L_0$-modules. By
Lemma 2.5, the decomposition $\pi_1(S_0) = \pi_1(K_0) +
\pi_1(L_0)$ can be considered in the matrix form (1).

Let $\varphi$ is a automorphism of $gl(2k)$, such that
$\varphi(X)=QXQ^{-1}$, where  $$Q= \frac{1}{\sqrt{2}}
\left(\begin{array}
                   { c  |  c    }
                     I_k  & I_k   \\ \hline
                     iI_k  & -iI_k    \\
           \end{array}
   \right)$$ where $I_k$ is an unit matrix of order $k$.

We extend $\varphi$ to an automorphism $\bar{\varphi}$ of
$gl(2k,2n)$ by the following formula:

$$\bar{\varphi} (Y) = \bar{Q}Y\bar{Q}^{-1}   \eqno(5) $$

 where  $$Y= \left(\begin{array}
                   { c  |  c    }
                     A  & B     \\ \hline
                     C  & D     \\
           \end{array}
   \right),$$  $$\bar{Q} = \left(\begin{array}
                   { c  |  c    }
                     Q  & 0     \\ \hline
                     0  & I_{2n}     \\
           \end{array}
   \right)$$

Therefore we can obtain a new decomposition $\bar{\varphi}(S) =
\bar{\varphi}(K) + \bar{\varphi}(L)$. Let $S' = \bar{\varphi}(S)$,
$K' = \bar{\varphi}(K)$ and $L' = \bar{\varphi}(L)$.

First we consider $K'$. Since $K'_0 = \bar{\varphi}(K_0)$ and
$\pi(K_0)$ have the form (3), the formula (5) gives us the
following matrix form of $K'_0$:
$$\left\{    \left(\begin{array}
                   { c  |  c    }
                     A  & 0     \\ \hline
                     0  & D     \\
           \end{array}
   \right) \right\}   \eqno(6)  $$  where $A$ is a  set of skewsymmetric  matrices with the
   first column and row  zero.

Next we want to show that $K'_1$ has the form
$$\left\{    \left(\begin{array}
                   { c  |  c    }
                     0  & B     \\ \hline
                     C  & 0     \\
           \end{array}
   \right) \right\} \eqno(7)$$
where the first  column of $C$ is zero.

Let $I_1 \cong o(2k-1)$, $I_2 \cong sp(2n)$ be ideals of $K'_0$.
By Lemma 2.2(d), $K'_1= [I_1,K'_1]$. Notice that $I_1$ has a form
(6) where $D=0$. Hence $[I_1,K'_1]$ takes the form
$$\left\{    \left(\begin{array}
                   { c  |  c    }
                     0  & AB     \\ \hline
                     -CA & 0     \\
           \end{array}
   \right) \right\} .$$

Since the first column of $A$ is zero, the first column of $CA$ is
zero. On the other hand $[I_1,K'_1]$ coincides with $K'_1$. Hence
we have proved that $K'_1$ has a form (7) where the first column
of $C$ is zero.

Next we consider $S'$. Let $\pi$ be a projection of  $(V \otimes
W^{\ast}) \oplus (V^{\ast} \otimes W)$ (or equivalently
$gl(2k,2n)_1$) onto $V^{\ast} \otimes W$. Since $S'_1$ is
$S'_0$-module, $\pi(S'_1)$ is also $S'_0$-module. We are going to
prove that $\pi(S'_1)$ coincides with $V^{\ast} \otimes W$.

First, $S'_0$-module $\pi(S'_1)$ is not zero. Indeed, if
$\pi(S'_1) = \{0\}$
 than $S'_1$ has the following matrix form:
$$\left\{     \left(\begin{array}
                   { c  |  c    }
                     0  & *     \\ \hline
                     0  & 0     \\
           \end{array}
   \right) \right\} .$$ Hence $[S'_1,S'_1] = \{0\}$. This contradicts the fact
   that, by Lemma 2.2(c), $[S'_1,S'_1] = S'_0 \neq \{0\}$.

   By Lemma 2.2(b), $S'_0$-module $S'_1$ is irreducible. Therefore
$S'_0$-module $S'_1$ is isomorphic to $S'_0$-module $\pi(S'_1)
\subseteq V^{\ast} \otimes W$. Since  dimension of $S'_1$ is
$4kn$, it follows that  dimension of
 $\pi(S'_1)$ is also $4kn$.  On the other hand,   $\dim({ V^{\ast} \otimes W}) =  (\dim{V})(\dim{W})=
4kn$. Hence $\pi(S'_1)$ coincides with $V^{\ast} \otimes W$. This
implies that $S'$ has the form (7), where $C$ is an arbitrary
matrix of order $2n \times 2k$.

Finally we consider  $L'= \bar{\varphi}(L)$. Let us assume the
contrary, that is, there are $\lambda \in F$ and $j_0$,  $j_1 \in
\{1 \ldots d\}$, $j_1\neq j_2$, such that $\varrho_{ij_0}(L_1)=
\lambda \varrho_{ij_1}(L_1)$ for any $x \in L_1$, $i \in \{1,2\}$.
The matrix realization of $L_1$ has the following form:

$$ \left\{
\left(\begin{array}
                   {          c        c      |       c       c      c    }
                                   &          &           &       &         \\
                          {\mbox{\Large 0}}&  &           &  {\Large \ast}    &         \\
                                   &          &           &       &         \\ \hline
                     M_{11}           &   M_{21}       &           &       &         \\
                      \vdots               &    \vdots      &           &  {\mbox{\Large 0}}     &         \\
                     M_{1d}          &   M_{2d}       &           &       &         \\

           \end{array}
   \right) \right\}
$$ where $M_{ij_0} = \lambda M_{ij_1} $ for any $i \in \{1,2\}$.

Using the formula (5), we obtain that  $L'_1$ has the following
form:
$$ \left\{
\left(\begin{array}
                   {          c        c      |       c       c      c    }
                                   &          &           &       &         \\
                          {\mbox{\Large 0}}&  &           &  {\Large \ast}    &         \\
                                   &          &           &       &         \\ \hline
                     M'_{11}         &   M'_{21}       &           &       &         \\
                      \vdots         &    \vdots      &           &  {\mbox{\Large 0}}     &         \\
                     M'_{1d}         &   M'_{2d}       &           &       &         \\

           \end{array}
   \right) \right\}
$$
where $M'_{ij} = M_{ij}Q^{-1}$. Hence $M'_{ij_0} = \lambda
M'_{ij_1} $ for any $i \in \{1,2\}$. On the other hand, $S'=K'+
L'$ and $\pi(S')=\pi(K')+ \pi(L')$. Since the first columns of
matrices from $\pi(S')$ are  arbitrary vectors from $F^{2n}$ and
the first columns of matrices from $\pi(K')$ are zero, it follows
that the first columns of matrices from $\pi(L')$ are arbitrary
vectors from $F^{2n}$. This contradicts the fact that $M'_{1j_0} =
\lambda M'_{1j_1}. $ \hfill $\fbox{}$

\end{proof}

\begin{lemma}
Let $S=K+L $ where $S\cong osp(2k,2n)$, $K\cong osp(2k-1,2n)$,
$L\cong sl(k,l)$. Then $L_0$-module $W_{j_0}$, $j_0 \in \{1 \ldots
d\}$ is neither of the type 1 no type 2.

\end{lemma}

\begin{proof}

We fix a basis in $V \oplus W$ of elements of subspaces $V_i$,
$i=1, 2$ and $W_j$, $j=1 \ldots d$. In this basis $L_0$ takes a
form
$$ \left\{
\left(\begin{array}
                   { c       c     |    c       c      c   }

                     Y    &   0    &        &       &          \\
                     0    & -Y^{t} &        &       &          \\ \hline

                          &         &  Z_1  &         &   {\mbox{\Large 0}}  \\
                          &         &        & \ddots  &        \\
                          &         & {\mbox{\Large 0}}   &         &  Z_d    \\
           \end{array}
   \right) \right\}  \eqno(8)
$$ where $Z_j$ is a matrix realization of $L_0$-module $W_j$

Let us assume the contrary, $L_0$-module $W_{j_0}$ is either of
the type 1 or type 2 for some $j_0$. Hence  $I_2$ acts trivially
on $W_{j_0}$. By Corollary 2.6, $I_2$ acts trivially on $V$.
Therefore $I_2$ has the form (8) where $Z_{j_0}=0$, $Y=0$.  Let
$L_1$ have the form:
$$ \left\{
\left(\begin{array}
                   { c              c      |       c       c      c    }
                                   &          &           &       &         \\
                          {\mbox{\Large 0}}&  &           &  {\Large \ast}    &         \\
                                   &          &           &       &         \\ \hline
                     M_{11}     &   M_{21}       &           &       &         \\
                      \vdots    &    \vdots      &           &  {\mbox{\Large 0}}     &         \\
                     M_{1d}     &   M_{2d}       &           &       &         \\

           \end{array}
   \right) \right\}.
$$
 Then $[I_2,L_1]$  has the form:
$$ \left\{
\left(\begin{array}
                   { c              c      |       c       c      c    }
                                   &          &           &       &         \\
                          {\mbox{\Large 0}}&  &           &  {\Large \ast}    &         \\
                                   &          &           &       &         \\ \hline
                     M'_{11}     &   M'_{21}       &           &       &         \\
                      \vdots    &    \vdots      &           &  {\mbox{\Large 0}}     &         \\
                     M'_{1d}     &   M'_{2d}       &           &       &         \\

           \end{array}
   \right) \right\}
$$  where $M'_{ij_{0}}$, $i \in
\{1,2\}$ are zero since $M'_{ij_{0}} =
Z_{j_0}M_{ij_0}-M_{ij_0}0=0$. On the other hand, by Lemma 2.1(d),
$[I_2,L_1]=L_1$. By Lemma 2.7, it contradict the fact that,
 there exists $i_0 \in \{1,2\}$ such that
$M_{i_0j_0}$ is not zero.  \hfill $\fbox{}$
\end{proof}

$\quad$

In this paper we will employ the following construction. Let
$L$-module $V(\lambda)$ and $L'$-module $V(\mu)$ are two
irreducible modules. Then one can define $L \oplus L'$-module
$V(\lambda) \otimes V(\mu)$ in the  natural way $$(X,Y)(v \otimes
w)= X(v) \otimes w + v \otimes Y (w). \eqno(9)$$ The following
theorem  (see \cite{GG}) holds:
\begin{lemma}
If $L$-module $V(\lambda)$ and $L'$-module $V(\mu)$ are two
irreducible modules then $L \oplus L'$-module $V(\lambda) \otimes
V(\mu)$ is also irreducible with a highest weight $(\lambda,\mu)$.
\end{lemma}

\begin{lemma}
Let $S=K+L$ where $S\cong osp(2k,2n)$, $K\cong osp(2k-1,2n)$,
$L\cong sl(k,l)$. If $L_0$-module $W_{j_0}$, $j_0 \in \{1 \ldots
d\}$ is of the type 3 then $I_2$-module $W_{j_0}$ is either
standard or dual.
\end{lemma}

\begin{proof}
By Lemma 2.7(a), there exists $i_{0}$ such that
$\pi_{i_{0}j_0}(L_1)\neq 0$. We consider $L_0$-module
$V_{i_0}^{\ast} \otimes W_{j_{0}}$. Since $I_1$-module $V_{i_0}$
and $I_2$-module $W_{j_{0}}$ are irreducible, by Lemma 2.9,
$L_0$-module $V_{i_0}^{\ast} \otimes W_{j_{0}}$ is also
irreducible. Therefore $L_0$-module $\pi_{i_{0}j_0}(L_1)$
coincides with $V_{i_0} \otimes W_{j_{0}}^{\ast}$. By Lemma
2.1(b), $L_0$-module $L_1$ is a direct sum of two irreducible
$L_0$-submodules of dimensions $kl$ each. Since
$\pi_{i_{0}j_0}(L_1)$ is irreducible $L_0$-module, it follows that
 dimension of $\pi_{i_{0}j_0}(L_1)$ is equal to $kl$ because. On
the other hand, we have  $$(\dim{V_{i_0}})(\dim{W_{j_0}}) = \dim{(
V_{i_0}^{\ast} \otimes W_{j_{0}})} = \dim{\pi_{i_{0}j_0}(L_1)} =
kl.$$

Since  $V_{i_0}$ is a nontrivial $sl(k)$-module and $W_{j_0}$ is a
nontrivial $sl(l)$-module it follows that dim $V_{i_0} \geq k$ and
dim $W_{j_0} \geq l$. Therefore dim $V_{i_0} = k$ and dim $W_{j_0}
= l$. Hence $W_{j_0}$ is either standard or dual. \hfill $\fbox{}$

\end{proof}

$\quad$

In the following Lemma  $L_0=I_1 \oplus I_2$ where $I_1$, $I_2$
are ideals of $L_0$.

\begin{lemma}
Let $U$ be an irreducible $L_0$-module such that $I_1(U)\neq 0 $
and $I_2(U)\neq 0 $. Then there exist $U'$, $U'' \subseteq U$ such
that $U'$ is an  irreducible $I_1$-module and $U''$ is an
irreducible $I_2$-module. Moreover, $U$ isomorphic to $U' \otimes
U''$ as $L_0$-module.
\end{lemma}

\begin{proof}

Let $\lambda = (\lambda',\lambda'')$ be the highest weight of
$L_0$-module $U$ where $\lambda'$ and $\lambda''$ correspond to
$I_1$ and $I_2$.  Next we can take  $I_1$-module $U_1$ and
$I_2$-module $U_2$ with the highest weights $\lambda'$ and
$\lambda''$, respectively.  We define $I_1 \oplus I_2$-module $U_1
\otimes U_2$ as was shown above (see (9)).  By Lemma 2.9, $I_1
\oplus I_2$-module $U_1 \otimes U_2$ is irreducible with the
highest weight $(\lambda', \lambda'') = \lambda$. Therefore $I_1
\oplus I_2$-modules $U_1 \otimes U_2$ and $U$ are isomorphic. Let
$\psi$ be a isomorphism between $U_1 \otimes U_2$ and $U$. Next we
choose $u_1\in U_1$ and $u_2\in U_2$. By formula (9), $U_1 \otimes
u_2$ is $I_1$-module and $u_1 \otimes U_2$ is $I_2$-module.
Moreover, $U_1 \otimes u_2$  isomorphic to $U_1$ as $I_1$-module
and $u_1 \otimes U_2$  isomorphic to $U_2$ as $I_2$-module. We
define $U'= \psi(U_1 \otimes u_2)$ and $U''= \psi(u_1 \otimes
U_2)$. Since $U_1 \cong U' $ as $I_1$-module and $U_2 \cong U'' $
as $I_2$-module, it follows that $U_1 \otimes U_2 \cong U' \otimes
U''$ as $I_1 \oplus I_2$-module. Therefore  $U$ isomorphic to $U'
\otimes U''$ as $L_0$-module. \hfill $\fbox{}$
\end{proof}

\begin{lemma}
Let $S=K+L $ where $S\cong osp(2k,2n)$, $K\cong osp(2k-1,2n)$,
$L\cong sl(k,l)$. Then for any $j_0 \in \{1 \ldots d\}$,
$L_0$-module $W_{j_0}$  is not of the type 4 .

\end{lemma}

\begin{proof}

 Let us assume the
contrary, that is, there exist $j_0$ such that $L_0$-module
$W_{j_0}$  is  of the type 4. By  Lemma 2.11, there exist
subspaces $W'_{j_0} \subseteq W_{j_0} $ and $W''_{j_0} \subseteq
W_{j_0}$ such that $W'_{j_0}$ is irreducible $I_1$-module,
$W''_{j_0}$ is irreducible $I_2$-module and $W_{j_0} \cong
W'_{j_0} \otimes W''_{j_0}$ as $I_1 \oplus I_2$-modules.

First we show that  dim $W'_{j_0} = k$ and dim $W''_{j_0} =l$.
Since $W'_{j_0}$ is an irreducible $sl(k)$-module and $W''_{i_0}$
is an irreducible $sl(l)$-module, it follows that  $\dim{W'_{j_0}}
\geq k$ and $\dim{W''_{j_0}} \geq l$, respectively. Without any
loss of generality, we assume that dim $W'_{j_0} > k$. Therefore $
2n = \dim{W} \geq \dim{W_{j_0}} = \dim{W'_{j_0}} \dim{W''_{j_0}} >
kl$. Since  $\dim{ S_1} \leq \dim{K_1} + \dim{L_1}$, it follows
that $ \dim{L_1} \geq \dim{ S_1} - \dim{K_1} = 2(2k)(2n) -
2(2k-1)(2n) = 4n $. On the other hand, $ \dim{L_1} = 2kl$. Hence
$2kl  \geq 4n$. This contradicts the fact that $2n > kl$.
Therefore $\dim{W'_{j_0}} = k$, $\dim{W''_{j_0}} = l$ and
$W=W_{j_0}$. If we denote $W'_{j_0}$ and $W''_{j_0}$ as $W'$ and
$W''$, then $W \cong W' \otimes W''$.

  Let us fix the following basis for $W$: $\{e'_{i} \otimes e''_{j} \}$, where
$\{e'_{i} \}$ is a basis of $W'$ and $\{e''_{j} \}$ is a basis of
$W''$. If we consider $W$ as $I_1$-module then it can be expressed
as the direct sum of $I_1$-modules $W' \otimes e''_k$:
$$ W = (W' \otimes e''_1) \oplus \ldots  \oplus  (W' \otimes e''_k).  \eqno(10) $$

Let us prove that the projection $\pi$ of $L_1$ onto $V^{\ast}
\otimes W$ is not zero. Indeed, if $\pi(L_1) = \{0\}$  then $L_1$
has the following matrix form:
$$\left\{     \left(\begin{array}
                   { c  |  c    }
                     0  & *     \\ \hline
                     0  & 0     \\
           \end{array}
   \right) \right\} .$$ Hence $[L_1,L_1] = \{0\}$. This contradicts the fact
   that, by Lemma 2.1(c), $[L_1,L_1] = L_0 \neq \{0\}$. Hence $\pi(L_1) \neq
   \{0\}$.  Therefore there exists $i_0 \in \{1, 2\}$ such
that the projection of $L_1$ onto $V_{i_0}^{\ast} \otimes W$ is
not zero. Let us consider $V_{i_0}^{\ast} \otimes W$ as
$I_1$-module. From (10) we obtain that
$$V_{i_0}^{\ast}\otimes W = (V_{i_0}^{\ast}\otimes(W' \otimes e''_1) ) \oplus \ldots \oplus
(V_{i_0}^{\ast}\otimes(W' \otimes e''_m) )  $$ where
$V_{i_0}^{\ast}\otimes(W' \otimes e''_k) $ are also $I_1$-modules.
The projection of $L_1$ onto $V_{i_0}^{\ast}\otimes(W' \otimes
e''_k) $  is not zero for some $k_0$ since the projection of $L_1$
onto $V_{i_0}^{\ast} \otimes W$ is not zero.

We consider  $I_1$-module $V_{i_0}^{\ast}\otimes(W' \otimes
e''_{k_0}) $. By Corollary 2.6, $I_1$-module $V_{i_0}$ is either
standard or dual. Since   $\dim{W'} = k$, it follows that
$I_1$-module $W'$ is also either standard or dual. Next we apply
 Young tableaux technique (see \cite{EL}) to find
irreducible submodules of $I_1$-module $(V_{i_0}^{\ast}\otimes W')
\otimes e''_{k_0} $.

 Let  $\varrho$ be either a standard or duel representation and
$\varrho'$ be  also either a standard or duel representation of
$sl(k)$. Then the tensor product $\varrho \otimes \varrho'$ is
also a representation of $sl(k)$ and, by Young tableaux technique,
it  can only contain  irreducible subrepresentations with highest
weights $(2,0, \ldots ,  0 )$, $(0,1, 0, \ldots , 0 )$, $(1,0,
\ldots , 0, 1 )$ or a trivial representation.

 Since $I_1$-modules $V_{i_0}$ and $W'$ are either standard or dual, we obtain that $I_1$-module
$(V_{i_0}^{\ast}\otimes W') \otimes e''_{k_0} $ can only contain
irreducible submodules with  highest weights listed above. On the
other hand, by Lemma 2.1(e) $I_1$-module $L_1$ has only standard
irreducible submodules of  dimension $k$. Contradiction.

\hfill $\fbox{}$

\end{proof}

\begin{lemma}
Let $S=K+L $ where $S\cong osp(2k,2n)$, $K\cong osp(2k-1,2n)$,
$L\cong sl(k,l)$. Then for any pairwise different $j_1,j_2 \in \{1
\ldots d\}$, $L_0$-module $W_{j_1}$ is not isomorphic to
$L_0$-module $W_{j_2}$.

\end{lemma}
\begin{proof}

Let us assume the contrary, that is, $L_0$-modules $W_{j_1}$ and
$W_{j_2}$ are isomorphic.

By Lemmas 2.8 and 2.12, any $L_0$-module $W_{j_1}$ is of the type
3. Moreover, by Lemma 2.10, $L_0$-module $W_{j_1}$ is either
standard or dual.

Without any loss of generality, we only consider the case than
$L_0$-module $W_{j_1}$ is standard. Hence $L_0$-module $W_{j_2}$
is also standard. By Corollary 2.6, $L_0$-module $V_1$ is standard
and $L_0$-module $V_2$ is dual. Therefore we obtain two cases:

1. $L_0$-modules $V_1^{\ast} \otimes W_{j_1}$ and  $V_1^{\ast}
\otimes W_{j_2}$ have the same highest weight $(\mu,\lambda)$,
where $\mu =(0, \ldots, 0,1 )$ and $\lambda =(1,0, \ldots, 0 )$

2. $L_0$-modules $V_2^{\ast} \otimes W_{j_1}$ and $V_2^{\ast}
\otimes W_{j_2}$ have the same highest weight $(\lambda,\lambda
)$, where $\lambda =(1,0, \ldots, 0 )$

By Lemma 2.1(b), $L_0$-module $L_1$ is a direct product of two
irreducible submodules with highest weights $(\mu,\lambda)$ and
$(\lambda,\mu)$ respectively. Hence the projections of $L_1$ onto
$V_2^{\ast} \otimes W_{j_1}$ and $V_2^{\ast} \otimes W_{j_2}$ are
zero since $L_0$-module $L_1$ contains no submodules with the
highest weights $(\lambda,\lambda ).$

Next, $L_0$-module $V_1^{\ast} \otimes W_{j_1}$  and $L_0$-module
$V_1^{\ast} \otimes W_{j_2}$  are irreducible and have the same
highest weights. Hence they are isomorphic as $L_0$-modules.
Therefore $\varrho_{ij_1}(V^{\ast} \otimes W)$ and
$\varrho_{ij_2}(V^{\ast} \otimes W)$ are also isomorphic as
$L_0$-modules. By Schur's Lemma the only endomorphisms between
these $L_0$-modules are scalars. However this contradicts to Lemma
2.7(b). \hfill $\fbox{}$

\end{proof}

\begin{theorem}
Let  $ S = {osp}(m,2n)$ be decomposed into the sum of two proper
simple subalgebras ${ K}$ and ${ L}$ of the type ${osp}(p,2q)$ and
${sl}(s,l)$, respectively. Then $m$ is even, $m=2k$ and $K \cong
osp(2k-1,2n)$, $L \cong sl(k,n)$.

\end{theorem}
\begin{proof}
By Lemma 2.5, $m=2k$ and $K \cong osp(2k-1,2n)$, $L \cong
sl(k,l)$. We only have to prove that $l=n$. Let us consider
$L_0$-modules $W = W_1 \oplus \ldots \oplus W_d$. By Lemmas 2.8
and 2.12, for any $j \in \{1 \ldots d\}$ $L_0$-module $W_j$ is not
of the type 1, 2 and 4. Hence any $L_0$-module $W_j$ is of the
type 3. Moreover, by Lemma 2.10, $I_2$-module $W_j$ has dimension
$l$. Since $\pi_2(I_2)\neq {0}$, $\pi_2(I_2) \subseteq sp(2n)$ and
$I_2 \cong sl(l)$, it follows that $l< 2n$. Hence $\dim{W_j}
=l<2n= \dim{W}$. Therefore $W$ contains at least two $L_0$-modules
$W_1$ and $W_2$ of type 3.

Next we show that $d=2$. Let us assume the contrary. There exists
$L_0$-modules $W_3$.  Since $L_0$-modules $W_3$ is of the type 3,
it follows that $L_0$-module $W_3$ is either standard or dual. By
Lemma 2.13, $L_0$-modules $W_1$ and $W_2$ are not isomorphic.
Therefore $L_0$-module $W_3$ is isomorphic to either $L_0$-modules
$W_1$ or $L_0$-modules $W_2$. However, this contradicts  Lemma
2.13. Since $d=2$ it follows that $l = \dim{W_1} = \dim{W}/2 = n
$. Therefore $L \cong sl(k,n)$. \hfill $\fbox{}$

\end{proof}

\begin{fact}

We consider  Lie superalgebra $S \cong osp(2k,2n)$ in the standard
matrix realization:
 $$\left\{ \left(\begin{array}
                   { c  |  c    }
                     A  & B     \\ \hline
                     C  & D     \\
           \end{array}
   \right) \right\}$$ where $A \in gl(2k)$ and $D \in gl(2n)$ and
   $A^{t}=-A$, $D^{t}G=-GD$, $C=GB^{t}$,  $G$ is given by $$G=
 \left(\begin{array}
                   { c  |  c    }
                     0     & I_n     \\ \hline
                     -I_n  & 0     \\
           \end{array}
   \right)$$

 Let the first subalgebra $K \cong osp(2k-1,2n)$ have the form:
$$
\left\{ \left(\begin{array}{c|ccc}
             0 &   0         & 0  & 0 \\ \hline
             0 &    &    &  \\
             0 &             & {\Huge X}   &  \\
             0 &             &    &  \\
             \end{array}\right) \right\}
$$ where $X$  is any $(2k+2n-1) \times (2k+2n-1)$ orthosymplectic  matrices.

The second subalgebra $L  \cong sl(k,n) $  consists of all
matrices of the form:
$$
\left\{ \left(\begin{array}{cc |cc}
                          E      &      -F        &      P         &    Q^t         \\
                          F     &     E           &    iP          &   -iQ^{t}      \\ \hline
                          Q     &      -iQ        &    D           &    0         \\
                         -P^t    &     -iP^{t}    &    0           &    -D^t        \\
             \end{array}\right) \right\} \eqno(11)
$$
where  $E$ is a skewsymmetric matrix  of order $k$, $F$ is a
symmetric matrix  of order $k$,
 $P$ is a  matrix  of order $k \times n$,
 $Q$ is a matrix  of order $n \times k$
 and $D$ is a  matrix  of order $n$ with zero trace.

 Then   $S=K+L$ is a decomposition of a simple Lie
superalgebra onto the sum of two simple subalgebras.

\end{fact}

\begin{proof}

First we prove that the set of matrices (11) forms $sl(k,n)$. The
standard matrix realization of $sl(k,n)$ has the form: $$  \left\{
\left(\begin{array}
                   { c  |  c    }
                     X  & P     \\ \hline
                     Q  & Y     \\
           \end{array}
   \right) \right\}$$ where   $X$ is a  matrix  of order $k$ with zero trace,
 $P$ is a  matrix  of order $k \times n$,
 $Q$ is a  matrix  of order $n \times k$,
   $Y$ is a  matrix  of order $n$ with zero trace. Then $sl(k,n)$
also has the following matrix realization:
   $$ \left\{  \left(\begin{array}
                   { c  |  c    }
                     -X^t  & Q^{t}     \\ \hline
                     -P^{t}  & -Y^t     \\
           \end{array}
   \right) \right\}$$

   Therefore we consider $L'\cong sl(k,n)$ in the  form:
$$ \left\{  \left(\begin{array}
                   { c    c      |  c       c    }
                     X  & 0      &  P   &   0    \\
                     0  & -X^{t} &  0   &   Q^t    \\ \hline
                     Q  & 0      &  Y   &    0     \\
                     0  & -P^t    &  0   &   -Y^t    \\
           \end{array}
   \right) \right\}$$

   Let $\bar{\varphi}$ be an automorphism of $gl(2k,2n)$ of the
   form (5). The direct calculation gives us that
   $\bar{\varphi}(L')$ has the form (11) where $E=A-A^{t}$,
   $F=i(A+A^t)$. Therefore the set of matrices of the form (11)
   forms  $sl(k,n)$.

   Next we prove that  the sum of two vector spaces $K$ and $L$
   coincides with $S$.

Let $$B=  \left(\begin{array}
                   { c  |  c    }
                     B_{11}  &  B_{12}     \\ \hline
                     B_{21}  &  B_{22}     \\
           \end{array}
   \right).$$

   Then
$$ C=GB^{t}= \left(\begin{array}
                   { c  |  c    }
                     B_{12}^{t}  &  B_{22}^{t}     \\ \hline
                     -B_{11}^{t}  &  -B_{21}^{t}     \\
           \end{array}  \right) $$
We set $B_{11}=P$ and $B_{12}=Q^{t}$. Then $B_{12}^{t} = Q$ and
$-B_{11}^{t}= -P^t$. Since $P$ and $Q$ are arbitrary matrices of
order $k \times n$ and $n \times k$, respectively, it follows that
the first raw and column of matrices from $L$
   coincides with the first raw and column of matrices from $S$.    \hfill $\fbox{}$

\end{proof}

$\quad$

\subsection{ Decompositions of  $osp(2k,2n)$ as
the sum of   $sl(p,q)$ and ${sl}(s,l)$}

In this section  we consider the decomposition of $osp(m,2n)$ into
the sum of two proper simple subalgebras ${ K} \cong sl(p,q)$ and
${ L} \cong {sl}(s,l) $.

\begin{theorem}
A Lie superalgebra  $ S \cong {osp}(m,2n)$ cannot  be decomposed
into the sum of two proper simple subalgebras ${ K}$ and ${ L}$ of
the  type $sl(p,q)$ and ${sl}(s,l)$, respectively.

\end{theorem}

\begin{proof}
By Lemma 2.2(a),  ${S}_0 = {o(m)} \oplus sp(2n)$. We define two
projections $\pi_1$ and $\pi_2$ of ${S}_0  $ onto the ideals
$o(m)$ and $sp(2n)$, $\pi_1:{ S}_0\to{o(m)}$ and $\pi_2:{ S}_0\to
{sp(2n)}$. Since  ${ K} \cong {sl}(p,q)$ and ${ L} \cong
{sl}(s,l)$, it follows that  ${ K}_0 \cong {sl(p)} \oplus {sl(q)}
\oplus { U}$ and  ${L}_0 \cong {sl(s)} \oplus {sl(l)} \oplus {
U}$. Since ${ K}_0$ and ${L}_0$ are reductive subalgebras, the
projections $\pi_1({ K}_0)$, $\pi_1({ L}_0)$, $\pi_2({ K}_0)$ and
$\pi_2({ L}_0)$ are also reductive as homomorphic images of
reductive algebras.

Since ${ S} = { K} + { L}$,  ${ S}_0$ is decomposable into the sum
of two subalgebras ${K}_0$ and ${ L}_0$, ${ S}_0 = { K}_0 + {
L}_0$. Therefore, $\pi_1({ S}_0) = \pi_1({ K}_0) + \pi_1({ L}_0)$
and $\pi_2({ S}_0) = \pi_2({ K}_0) + \pi_2({ L}_0)$, where
$\pi_1({\cal S}_0)={ o(m)}$ and $\pi_2({\cal S}_0)={ sp(2n)}$. We
have the decompositions of simple Lie algebras $o(m)$ and $sp(2n)$
into the sum of two reductive subalgebras.

By  Onichshik's Theorem, ${sp(2n)}$ and  $o(m)$  cannot be
decomposed into the sum of two subalgebras of this form. Therefore
$ S \cong {osp}(m,2n)$ cannot be decomposed into the sum of  ${ K}
\cong {sl}(p,q)$ and ${ L} \cong {sl}(s,l)$\hfill $\fbox{}$

\end{proof}

\subsection{ Properties of subalgebras ${osp}(p,2q)$ and ${osp}(s,2l)$ in
the decomposition }

We consider the decomposition of $osp(m,2n)$ into the sum of two
proper simple subalgebras ${ K} \cong {osp}(p,2q)$ and ${ L} \cong
{osp}(s,2l) $.

\begin{lemma}

 Let  ${ S}$ be a Lie superalgebra of the type
${osp}(m,2n)$, and $S$ be decomposed into the sum of two proper
simple subalgebras ${ K}$ and ${ L}$ of the  type ${osp}(p,2q)$
and ${osp}(s,2l)$, respectively. Then two cases are possible:

1.  $m=4k$  and  ${ K} \cong osp(4k-1,2n)$, ${ L}\cong osp(s,2k)$

2.  ${ K} \cong osp(p,2n)$, ${ L}\cong osp(m,2l)$ or ${ K} \cong
osp(m,2q)$, ${ L}\cong osp(s,2n)$

\end{lemma}

\begin{proof}

By Lemma 2.2(a),  ${S}_0 = {o(m)} \oplus sp(2n)$. We define two
projections $\pi_1$ and $\pi_2$ of ${S}_0  $ onto the ideals
$o(m)$ and $sp(2n)$ as follows $\pi_1:{ S}_0\to{o(m)}$ and
$\pi_2:{ S}_0\to {sp(2n)}$.  Since  ${ K} \cong {osp}(p,2q)$ and
${ L} \cong {osp}(k,2l)$, it follows that ${ K}_0 \cong {o(p)}
\oplus {sp(2q)}$ and ${ L}_0 \cong {o(k)} \oplus {sp(2l)}$. Both
 ${ K}_0$ and ${L}_0$ are semisimple. Hence
 $\pi_1({ K}_0)$, $\pi_1({ L}_0)$, $\pi_2({ K}_0)$ and
$\pi_2({ L}_0)$ are also semisimple as homomorphic images of
semisimple algebras.

Since ${ S} = { K} + { L}$,  ${ S}_0$ is decomposable into the sum
of two subalgebras ${K}_0$ and ${ L}_0$, ${ S}_0 = { K}_0 + {
L}_0$. Therefore, $\pi_1({ S}_0) = \pi_1({ K}_0) + \pi_1({ L}_0)$
and $\pi_2({ S}_0) = \pi_2({ K}_0) + \pi_2({ L}_0)$. Moreover,
$\pi_1({\cal S}_0)={ o(m)}$ and $\pi_2({\cal S}_0)={ sp(2n)}$. Now
we have the decompositions of simple Lie algebras $o(m)$ and
$sp(2n)$ into the sum of two semisimple  subalgebras.

By Onichshik's Theorem, ${sp(2n)}$ has no decompositions  into the
sum of two proper reductive subalgebras of these types. Hence
$sp(2n) = \pi_2({ K}_0) + \pi_2({ L}_0)$ is a  trivial
decomposition and either $\pi_2({ K}_0)= sp(2n)$ or $\pi_2({
L}_0)= sp(2n)$. Without any loss of generality, we assume that
$\pi_2({ K}_0)= sp(2n)$. Hence $q=n$.

By Onichshik's Theorem,  $o(m)$ has two decompositions into the
sum of two proper reductive subalgebras:

1. If $m=2k$ then  $o(2k)=o(2k-1)+ sl(k),$

2. If $m=4k$ then $o(4k)=o(4k-1)+sp(2k).$

 The decomposition $o(m) = \pi_1({ K}_0) + \pi_1({ L}_0)$ cannot
 be of the first type, because $\pi_1({ K}_0)$ and $\pi_1({ L}_0)$ are not isomorphic to
 $sl(k)$.

Next the two cases occur:

1. The decomposition $o(m) = \pi_1({ K}_0) + \pi_1({ L}_0)$ has
the second form.

2. The decomposition $o(m) = \pi_1({ K}_0) + \pi_1({ L}_0)$ is
trivial.

In the first case either $\pi_1({ K}_0) \cong o(4k-1)$ or $\pi_1({
K}_0) \cong sp(2k)$. Let $\pi_1({ K}_0) \cong sp(2k)$. We have
that $\pi_2({ K}_0)= sp(2n)$. Hence ${ K}_0 \cong {sp(2k)} \oplus
{sp(2n)}$ or  ${ K}_0 \cong {sp(2n)}$.  This contradicts the fact
that  ${ K}_0 \cong {o(p)} \oplus {sp(2q)}$. Therefore $\pi_1({
K}_0) \cong o(4k-1)$ and $\pi_1({ L}_0) \cong sp(2k)$.  Since ${
K}_0 \cong {o(p)} \oplus {sp(2q)}$ and ${ L}_0 \cong {o(s)} \oplus
{sp(2l)}$ it follows that $p=4k-1$ and $l=k$.

In the second case either $\pi_1({ K}_0)= o(m)$ or $\pi_1({ L}_0)=
o(m)$. Let $\pi_1({ K}_0)= o(m)$. Since  $\pi_2({ K}_0)= sp(2n)$
it follows that $K_0$ coincides  with $S_0$. This contradicts the
fact that $K$ is proper  subalgebra of $S$. Therefore $\pi_1({
L}_0)= o(m)$.  Since ${ L}_0 \cong {o(s)} \oplus {sp(2l)}$ it
follows that $s=m$. \hfill $\fbox{}$

\end{proof}

 $\quad$

\begin{corollary}
Let $S=K+L$ and ${ K} \cong osp(4k-1,2n)$, ${ L}\cong osp(s,2k)$.
Then $L_0$-module $V$ is a direct sum of two $L_0$-modules of type
2, $V=V_1 \oplus V_2$. Moreover $I_2$ acts trivially on $V$, and
$I_1$-modules $V_1$, $V_2$ are standard.
\end{corollary}

\begin{proof}
From the previous Lemma we know that in the decomposition $o(4k) =
\pi_1({ K}_0) + \pi_1({ L}_0)$ the first component isomorphic to
$o(4k-1)$ and the second one isomorphic to $sp(2k)$. By Lemma 2.4,
there exists a basis of $V$ such that the decomposition
$o(4k)=o(4k-1)+ sp(2k)$ takes the matrix form (4). Hence $\pi_1({
L}_0)$ takes the form:
$$ \left\{    \left(\begin{array}
                   { c  |  c    }
                     Y  & 0     \\ \hline
                     0  & -Y^{t}     \\
           \end{array}.
   \right) \right\} $$ where $Y \in sp(2k)$.

Since $I_2$ acts trivially on $V$,  we obtain  that $V=V_1 \oplus
V_2$ and $L_0$-modules $V_1$, $V_2$ are of type 2
 \hfill $\fbox{}$

\end{proof}

\begin{corollary}
Let $S=K+L$ and  ${ K} \cong osp(p,2n)$, ${ L}\cong osp(m,2l)$.
Then $L_0$-module $V$ is irreducible of type 2. Moreover $I_2$
acts trivially on $V$, and $I_1$-modules $V$ is standard.
\end{corollary}

\begin{proof}
The proof follows from the fact that $\pi_1({ L}_0)$ coincides
with $\pi_1({ S}_0).$
\end{proof}

\subsection{ Decompositions of  $osp(4k,2n)$ as
the sum of \\  $osp(4k-1,2q)$ and $osp(s,2k)$}

In this section we  consider the case when $m=4k$  and ${ K} \cong
osp(4k-1,2q)$, ${ L}\cong osp(s,2k)$. Let $L_0 = I_1 \oplus I_2$
where $I_1 \cong sp(2k)$ and $I_2 \cong o(s)$.

\begin{lemma}
Let $S=K+L $ where $S\cong osp(4k,2n)$, $K\cong osp(4k-1,2n)$,
$L\cong osp(s,2k)$.

\noindent (a) There is no  $j_{0} \in \{1 \ldots d\}$ such that
for any $i \in \{1,2\}$, $\pi_{ij_{0}}(L_1)=0$ .

\noindent (b) There are no $j_0,j_1 \in \{1 \ldots d\}$ such that
for some $\lambda \in F$, $\varrho_{ij_{0}}(x)=
\lambda\varrho_{ij_{1}}(x)$ where $x \in L_1$, $i \in \{1,2\}$.

\end{lemma}

\begin{proof}

 Let us assume the contrary, that is, there exist $j_0 $, $j_1 \in \{1 \ldots d\}$ such
that for any $i \in \{1 \ldots r\}$ and $\lambda \in$ F,
$\pi_{ij_0}(L_1)= \lambda\pi_{ij_1}(L_1)$.

We choose a basis in $V \oplus W$ from elements of subspaces
$V_i$, $i=1 \ldots r$ and $W_j$, $j=1 \ldots d$. Let us identify
the elements form $S$ with their matrix realizations in this
basis.

By Lemma 2.4, the  $\pi_1(S_0) = \pi_1(K_0) + \pi_1(L_0)$ can be
considered in the  matrix form (4). Let $\bar{\varphi}$ be the
automorphism of $gl(4k,2n)$ defined in Lemma 2.7. Then we obtain a
new decomposition $\bar{\varphi}(S) = \bar{\varphi}(K) +
\bar{\varphi}(L)$. Let $S' = \bar{\varphi}(S)$, $K' =
\bar{\varphi}(K)$ and $L' = \bar{\varphi}(L)$.

 Acting in the same
matter as in Lemma  2.7 we obtain that  $L'_1$ has the following
form:
$$
\left(\begin{array}
                   { c         c        c      |       c       c      c    }
                                   &          &           &       &         \\
                          {\mbox{\Large 0}}&  &           &  {\Large \ast}    &         \\
                                   &          &           &       &         \\ \hline
                     M'_{11}          &   M'_{21}       &           &       &         \\
                      \vdots          &    \vdots      &           &  {\mbox{\Large 0}}     &         \\
                     M'_{1d}          &   M'_{2d}       &           &       &         \\

           \end{array}
   \right),
$$ where $M'_{ij_0} = \lambda M'_{ij_1} $ for any $i \in \{1, 2\}$.

Since $S'=K'+ L'$, it follows that $\pi(S')=\pi(K')+ \pi(L')$. As
shown in Lemma 2.7,   the first columns of matrices from $\pi(S')$
are arbitrary vectors from $F^{2n}$, and the first columns of
matrices from $\pi(K')$ are zero. Hence the first columns of
matrices from $\pi(L')$ are arbitrary vectors from $F^{2n}$. This
contradicts the fact that $M'_{1j_0} = \lambda M'_{1j_1} $ \hfill
$\fbox{}$

\end{proof}

\begin{lemma}
Let $S=K+L $ where $S\cong osp(4k,2n)$, $K\cong osp(4k-1,2n)$,
$L\cong osp(s,2k)$. Then for any $j_0 \in \{1 \ldots d\}$,
$L_0$-module $W_{j_0}$  is neither of the type 1 nor type 2.

\end{lemma}

\begin{proof}

We choose a basis in $V \oplus W$ from elements of subspaces
$V_i$, $i=1,2$ and $W_j$, $j=1 \ldots d$. In this basis $L_0$
takes the form
$$ \left\{
\left(\begin{array}
                   { c       c     |    c       c      c   }

                     Y    &   0    &        &       &          \\
                     0    & -Y^{t} &        &       &          \\ \hline

                          &         &  Z_1  &         &   {\mbox{\Large 0}}  \\
                          &         &        & \ddots  &        \\
                          &         & {\mbox{\Large 0}}   &         &  Z_d    \\
           \end{array}
   \right) \right\} ,\eqno (12)
$$ where $Z_j$ is a matrix realization of $L_0$-module $W_j$

Let us assume the contrary, $L_0$-module $W_{j_0}$ is either of
the type 1 or type 2. Hence  $I_2$ acts trivially on $W_{j_0}$. By
Corollary 2.17, $I_2$ acts trivially on $V$. Therefore $I_2$ has
the form (12) where $Z_{j_0}=0$, $Y=0$.  Let $L_1$ have the form:
$$ \left\{
\left(\begin{array}
                   {          c        c      |       c       c      c    }
                                   &          &           &       &         \\
                          {\mbox{\Large 0}}&  &           &  {\Large \ast}    &         \\
                                   &          &           &       &         \\ \hline
                     M_{11}           &   M_{21}       &           &       &         \\
                      \vdots          &    \vdots      &           &  {\mbox{\Large 0}}     &         \\
                     M_{1d}           &   M_{2d}       &           &       &         \\

           \end{array}
   \right) \right\},
$$
 Then $[I_2,L_1]$  has the form:
$$ \left\{
\left(\begin{array}
                   {          c        c      |       c       c      c    }
                                   &          &           &       &         \\
                          {\mbox{\Large 0}}&  &           &  {\Large \ast}    &         \\
                                   &          &           &       &         \\ \hline
                     M'_{11}        &   M'_{21}       &           &       &         \\
                      \vdots              &    \vdots      &           &  {\mbox{\Large 0}}     &         \\
                     M'_{1d}           &   M'_{2d}       &           &       &         \\

           \end{array}
   \right)\right\} ,
$$  where $M'_{ij_{0}}$, $i \in \{1, 2\}$ is zero since
$M'_{ij_{0}} = Z_{j_0}M_{ij_0}-M_{ij_0}0=0$. On the other hand, by
Lemma 2.2(d), $[I_2,L_1]=L_1$.  This contradicts the fact that, by
Lemma 2.19(a),
 there exists $i_0 \in \{1, 2\}$ such that
$M_{i_0j_0}$ is not zero.  \hfill $\fbox{}$

\end{proof}

\begin{lemma}
Let $S=K+L $ where $S\cong osp(4k,2n)$, $K\cong osp(4k-1,2n)$,
$L\cong osp(s,2k)$. Then  for any $j \in \{1 \ldots d\}$,
$L_0$-module $W_{j_0}$ is not of the type 4.
\end{lemma}

\begin{proof}

 Let us assume the
contrary, that is, there exist $j_0$ such that $L_0$-module
$W_{j_0}$  is  of the type 4. By  Lemma 2.11, there exist
subspaces $W'_{j_0} \subseteq W_{j_0} $ and $W''_{j_0} \subseteq
W_{j_0}$ such that $W'_{j_0}$ is irreducible $I_1$-module,
$W''_{j_0}$ is irreducible $I_2$-module and $W_{j_0} \cong
W'_{j_0} \otimes W''_{j_0}$.

First we show that  dim $W'_{j_0} = 2k$ and dim $W''_{j_0} =s$.
Since $W'_{j_0}$ is an irreducible $sp(2k)$-module and $W''_{i_0}$
is an irreducible $o(s)$-module, it follows that $\dim{W'_{j_0}}
\geq 2k$ and $\dim{W''_{j_0}} \geq s$, respectively. Without any
loss of generality, we assume that dim $W'_{j_0} > 2k$. Hence $ 2n
= \dim{W} \geq \dim{W_{j_0}} = \dim{W'_{j_0}} \dim{W''_{j_0}} >
2ks$. Since  $\dim{ S_1} \leq \dim{K_1} + \dim{L_1}$, it follows
that $ \dim{L_1} \geq \dim{ S_1} - \dim{K_1} \geq 2nm - 2n(m-1) =
2n
> 2ks$. This contradicts the fact that  $ \dim{L_1} = 2ks$ since $L \cong
osp(s,2k)$. Therefore $\dim{W'_{j_0}} = 2k$, $\dim{W''_{j_0}} = s$
and $W=W_{j_0}$. If we denote $W'_{j_0}$ and $W''_{j_0}$ as $W'$
and $W''$, then $W \cong W' \otimes W''$.

Next we identify $W$ with $W' \otimes W''$.   Let us fix the
following basis for $W$: $\{e'_{i} \otimes e''_{j} \}$, where
$\{e'_{i} \}$ is a basis of $W'$ and $\{e''_{j} \}$ is a basis of
$W''$. If we consider $W$ as $I_1$-module then it can be expressed
as the direct sum of $I_1$-modules $W' \otimes e''_k$:
$$ W = (W' \otimes e''_1) \oplus \ldots  \oplus  (W' \otimes e''_k).  \eqno(13) $$

Clearly  the projection of $L_1$ onto $V^{\ast} \otimes W$ is not
zero. Therefore there exists $i_0 \in \{1,2\}$ such that the
projection of $L_1$ onto $V_{i_0}^{\ast} \otimes W$ is not zero.
Let us consider $V_{i_0}^{\ast} \otimes W$ as $I_1$-module. From
(13) we obtain that
$$V_{i_0}^{\ast}\otimes W = (V_{i_0}^{\ast}\otimes(W' \otimes e''_1) ) \oplus \ldots \oplus
(V_{i_0}^{\ast}\otimes(W' \otimes e''_m) )  $$ where
$V_{i_0}^{\ast}\otimes(W' \otimes e''_k) $ are also $I_1$-modules.
The projection of $L_1$ onto $V_{i_0}^{\ast}\otimes(W' \otimes
e''_k) $  is not zero for some $k_0$ since the projection of $L_1$
onto $V_{i_0}^{\ast} \otimes W$ is not zero.

We consider  $I_1$-module $V_{i_0}^{\ast}\otimes(W' \otimes
e''_{k_0}) $. By Corollary 2.17, $I_1$-module $V_{i_0}$ is
standard. We have already proved that $I_1$-module $W'$ is
standard with  highest weight $(1,0, \ldots , 0 )$. Next we apply
generalized Young tableaux technique (see \cite{EL}) to find
irreducible submodules of $I_1$-module $(V_{i_0}^{\ast}\otimes W')
\otimes e''_{k_0} $.

If $\varrho$ and  $\varrho'$ are standard representations of
$sp(2k)$ ($o(k)$) with the same highest weight $(1,0, \ldots , 0
)$ then the tensor product $\varrho \otimes \varrho'$ is also a
representation of $sp(2k)$ ($o(k)$). It can be decomposed into the
direct sum of irreducible representations: $$ \varrho \otimes
\varrho' = \varrho_1 \oplus \varrho_2 \oplus \varrho_3 \eqno(14)$$
where $\varrho_1$ has  highest weight $(2,0, \ldots ,  0 )$,
$\varrho_2$ has  highest weight $(0,1, 0, \ldots , 0 )$ and
$\varrho_3$ is a trivial representation.

Therefore $I_1$-module $(V_{i_0}^{\ast}\otimes W') \otimes
e''_{k_0} $ contains only submodules with highest weights $(2,0,
\ldots ,  0 )$ and $(0,1, 0, \ldots , 0 )$. This contradicts the
fact that, by Lemma 2.2(e), $I_1$-module $L_1$ has only standard
irreducible submodules of  dimension $2k$. \hfill $\fbox{}$

\end{proof}

\begin{lemma}

Let $S=K+L $ where $S\cong osp(4k,2n)$, $K\cong osp(4k-1,2n)$,
$L\cong osp(s,2k)$. If $L_0$-module $W_{j_0}$, $j_0 \in \{1 \ldots
d\}$ is of the type 3 then $I_2$-module $W_{j_0}$ is standard.
\end{lemma}

\begin{proof}
By Lemma 2.19, there exists $i_{0}$ such that
$\pi_{i_{0}j_0}(L_1)\neq 0$. We consider $L_0$-module
$V_{i_0}^{\ast} \otimes W_{j_{0}}$. By Lemma 2.9, since
$I_1$-module $V_{i_0}$ and $I_2$-module $W_{j_{0}}$ are both
irreducible, $L_0$-module $V_{i_0}^{\ast} \otimes W_{j_{0}}$ is
irreducible. Therefore $L_0$-module $\pi_{i_{0}j_0}(L_1)$
coincides with $V_{i_0} \otimes W_{j_{0}}^{\ast}$ since
$\pi_{i_{0}j_0}(L_1)\neq \{ 0 \}$. By Lemma 2.2(b), $L_0$-module
$L_1$ is  irreducible, $\dim{L_1}=2ks$. Since
$\pi_{i_{0}j_0}(L_1)$ is irreducible $L_0$-module, dimension of
$\pi_{i_{0}j_0}(L_1)$ is $2ks$. On the other hand, we have
$$(\dim{V_{i_0}})(\dim{W_{j_0}}) = \dim{( V_{i_0}^{\ast} \otimes
W_{j_{0}})} = \dim{\pi_{i_{0}j_0}(L_1)} = 2ks.$$

Since  $V_{i_0}$ is a nontrivial $sp(2k)$-module and $W_{j_0}$ is
a nontrivial $o(s)$-module,  dim $V_{i_0} \geq 2k$ and dim
$W_{j_0} \geq s$. Therefore dim $V_{i_0} = 2k$ and dim $W_{j_0} =
s$. Hence $I_2$-module $W_{j_0}$ is standard. \hfill $\fbox{}$

\end{proof}

$\quad$

\begin{lemma}
Let $S=K+L $ where $S\cong osp(4k,2n)$, $K\cong osp(4k-1,2n)$,
$L\cong osp(s,2k)$. Then for any $j \in \{1 \ldots d\}$,
$L_0$-module $W_{j}$  is not of the type 3.

\end{lemma}

\begin{proof}

We consider $L_0$-module $W_1 \oplus \ldots \oplus W_d$. By Lemmas
2.20 and 2.21, $L_0$-module $W_j$ is not of the type 1, 2 and 4.
Hence any $L_0$-module $W_j$ is of the type 3. By Lemma 2.22,
$L_0$-module $W_j$ has  dimension $s$. Therefore $\dim{W_{j}} <
\dim{W}$ since $s<2n$. It follows that $W$ contains at least two
$L_0$-modules $W_1$ and $W_2$ of type 3.

By Corollary 2.17, $V$ is a direct sum of two $I_1$-modules $V_1$
and $V_2$. We consider $I_1 \oplus I_2$-module $V^{\ast}_1 \otimes
W_1$. Without any loss of generality, $\pi_{11}(L_1)\neq 0.$

 At first we prove that there exist $\lambda \in F$ such that  $\varrho_{12}(x)=
\lambda\varrho_{11}(x)$ for any $x \in L_1$ (see section 2.2). If
$\pi_{12}(L_1)= \{0\}$ we choose $\lambda =0$. Let
$\pi_{12}(L_1)\neq 0$. Since $L_0$-module $L_1$ is irreducible,
$L_0$-module $L_1$ isomorphic to $L_0$-modules $\pi_{11}(L_1)$ and
$\pi_{12}(L_1)$. Hence   $\pi_{11}(L_1)$ and $\pi_{12}(L_1)$ are
isomorphic as $L_0$-modules. Therefore $\varrho_{11}(L_1)$ and
$\varrho_{12}(L_1)$ are also isomorphic as $L_0$-modules. By
Schur's Lemma the only endomorphisms between these $L_0$-modules
are scalars. Hence  $\varrho_{12}(x)= \lambda\varrho_{11}(x)$ for
any $x \in L_1$

Next we prove that $\varrho_{21}(L_1) = \{0\} $ and $
\varrho_{22}(L_1)=\{0 \}$.

Since $I_1$ acts trivially on $W_1$ and $W_2$,  $I_1$ has the
form:
$$\left\{ \left(\begin{array}
                   { c  |  c    }
                     A  & 0     \\ \hline
                     0  & 0     \\
           \end{array}
   \right) \right\}
$$ where
$$A = \left\{     \left(\begin{array}
                   { c  |  c    }
                     Y  & 0     \\ \hline
                     0  & -Y^{t}     \\
           \end{array}
   \right) \right\}
$$ where $Y$ is a set of matrices  defined in Lemma 2.4.

Let $\langle A \rangle$ be associative enveloping algebra
generated by all matrix from $A$. Since  $\{Y\}$ is an irreducible
set and $Y \neq -Y^{t}$ for some $Y$, it follows that $\langle A
\rangle$ takes a matrix form
$$A = \left\{     \left(\begin{array}
                   { c  |  c    }
                     Y'  & 0     \\ \hline
                     0  & Y''     \\
           \end{array}
   \right) \right\}
$$ where $Y'$, $Y''$ are arbitrary matrices of order $2k \times 2k.$

Hence $\langle A \rangle$ contains the following  matrix
  $$J = \left(\begin{array}
                   { c  |  c    }
                     I  & 0     \\ \hline
                     0  & 0     \\
           \end{array}
   \right)
$$

Since
$$
  \left[\left(\begin{array}
                   { c  |  c    }
                     A  & 0     \\ \hline
                     0  & 0     \\
           \end{array}
   \right),\left(\begin{array}
                   { c  |  c    }
                     0  & *     \\ \hline
                     C  & 0     \\
           \end{array}
   \right)\right] =
   \left(\begin{array}
                   { c  |  c    }
                     0  & *     \\ \hline
                     -CA  & 0     \\
           \end{array}
   \right),
$$ we obtain  that $L_1$ contains a set of matrices
$$\left\{
  \left(\begin{array}
                   { c  |  c    }
                     0  &  \ast    \\ \hline
                     CJ  & 0     \\
           \end{array}
   \right) \right\} \eqno(15)
$$ where
$$CJ=     \left(\begin{array}
                   { c    c    }
                     C_{11}  & 0     \\
                     \vdots  & \vdots     \\
                     C_{s1}  & 0    \\
           \end{array}
   \right)
$$

The subspace of $L_1$ of the form (15) has  dimension $2ks$ since
$C_{11}$ can be any matrix of order $s \times 2k$. On the other
hand, $\dim{L_1}= 2ks$. Hence $L_1$ has the form (15).

Since $\pi_{21}(L_1)= \{0\}$ and $\pi_{22}(L_1)=\{0\}$, it follows
that
 $\varrho_{21}(L_1) = \{0\} $ and  $\varrho_{22}(L_1) = \{0\} $.
 Therefore $\varrho_{12}(x) = \lambda \varrho_{11}(x)$ and $\varrho_{22}(x) =
\lambda \varrho_{21}(x)$ for any $x \in L_1$. This contradicts
Lemma 2.19. \hfill $\fbox{}$
\end{proof}

\begin{corollary}
 A Lie superalgebra  $ S \cong {osp}(2k,2n)$ cannot  be decomposed
into the sum of two proper simple subalgebras ${ K}$ and ${ L}$ of
the  types $osp(4k-1,2q)$ and $osp(s,2k)$, respectively.

\end{corollary}

\begin{proof}
The proof follows from the fact that for any $j \in \{1 \ldots
d\}$, $L_0$-module $W_{j}$ is not of the type 1, 2, 3 and 4.
\hfill $\fbox{}$
\end{proof}

$\quad$ $\quad$

\subsection{ Decompositions of  $osp(m,2n)$ as
the sum of   $osp(p,2n)$ and $osp(m,2l)$}

In this section  we  consider the case when $S=K+L$ where ${ S}
\cong {osp}(m,2n)$, ${ K} \cong osp(p,2n)$ and ${ L}\cong
osp(m,2l)$.  Let $L_0 = I_1 \oplus I_2$ where $I_1 \cong o(m)$ and
$I_2 \cong sp(2l)$.

In the following Lemma we show that one of two subalgebras $K$ and
$L$, for example $L$, does not contain $L_0$-module $W_{j}$ of the
type 1 and 2.

\begin{lemma}

Let $S=K+L$ where ${ S} \cong {osp}(m,2n)$, ${ K} \cong osp(p,2n)$
and ${ L}\cong osp(m,2l)$.  Then, without any loss of generality,
$L_0$-module $W_{j}$, $j \in \{1 \ldots d\}$ is neither  of the
type 1 nor 2.

\end{lemma}

\begin{proof}

We choose a basis in $V \oplus W$ from elements of subspaces $V$
and $W_j$, $j=1 \ldots d$. In this basis $L_0$ takes the form
$$ \left\{
\left(\begin{array}
                        {      c   |    c       c      c   }
                               Y   &        &         &          \\ \hline
                                   &  Z_1   &         &   {\mbox{\Large 0}}  \\
                                   &        & \ddots  &        \\
                  &      {\mbox{\Large 0}}   &         &  Z_d    \\
           \end{array}
   \right) \right\}  \eqno (16)
$$ where $Z_j$ is a matrix realization of $L_0$-module $W_j$

Let us assume the contrary, that is, there exists $j_0$ such that
$L_0$-module $W_{j_0}$ is either of the type 1 or 2. Without any
loss of generality, let $j_0=d$.

 Hence  $I_2$ acts trivially on $W_{d}$. By
Corollary (2.18), $I_2$ acts trivially on $V$. Therefore $I_2$ has
the form (16) where $Z_{d}=0$, $Y=0$.  Let $L_1$ have the form:
$$ \left\{
 \left(\begin{array}
                   { c  |  c            c               c }
                     0  &  N_{11}&   \ldots          &  N_{1d}    \\ \hline
               M_{11}   &        &                   &    \\
                \vdots  &        & {\mbox{\Large 0}} &   \\
                  M_{1d}&        &                   &   \\

           \end{array}
   \right) \right\}. \eqno (17)
$$

 Then $[I_2,L_1]$ also has the form (17), where $M_{1d}$ and $N_{1d}$ are
zero since $Z_{d}M_{1d}-M_{1d}0=0$ and $0N_{1d}-N_{1d}Z_d=0$. On
the other hand, by Lemma 2.2(d), $[I_2,L_1]=L_1$. Therefore
$M_{1d}$ and $N_{1d}$ are zero and $L$ consists of the matrices
with the last row and column are zero.

   Acting in the same manner as above, we obtain
   that  $K\cong osp(p,2n)$ consists of the matrices with the
first row and column are zero. This  contradicts the fact that
$S=K+L$. \hfill $\fbox{}$

\end{proof}

\begin{lemma}

Let $S=K+L$ where ${ S} \cong {osp}(m,2n)$, ${ K} \cong osp(p,2n)$
and ${ L}\cong osp(m,2l)$. Then  for any $j \in \{1 \ldots d\}$,
$L_0$-module $W_{j}$ is not of the type 4.
\end{lemma}
\begin{proof}

The proof of this Lemma is similar to the proof of Lemma 2.21.
\end{proof}

\begin{lemma}
Let $S=K+L$ where ${ S} \cong {osp}(m,2n)$, ${ K} \cong osp(p,2n)$
and ${ L}\cong osp(m,2l)$. If $L_0$-module $W_{j_0}$, $j_0 \in \{1
\ldots d\}$ is of the type 3 then $I_2$-module $W_{j_0}$ is
standard.
\end{lemma}

\begin{proof}

Without any loss of generality, we only consider the case $j_0=1$.

First we show that $\pi_{11}(L_1)\neq \{ 0\}$. Let us assume the
contrary, that is, $\pi_{11}(L_1) = \{0\}$.  Then $L_1$ takes the
form (17) where $M_{11}$ is zero. Hence $[L_1, L_1]$ takes the
form (16) where  $Z_1$ is zero.

 On the other hand, by Lemma 2.2(b), $L_0 =
[L_1, L_1]$. Hence $L_0$ also have the form (16) where $Z_{1}$ is
zero. This contradicts the fact that $L_0$-module $W_1$ is of the
type 3. Therefore $\pi_{11}(L_1)\neq \{ 0\}$.

Next we consider $L_0$-module $V_{1}^{\ast} \otimes W_{1}$. By
Lemma  2.9, $L_0$-module $V_{1}^{\ast} \otimes W_{1}$ is
irreducible since $I_1$-module $V_{1}$ and $I_2$-module $W_{1}$
are irreducible. Therefore $L_0$-module $\pi_{11}(L_1)$ coincides
with $V_{1} \otimes W_{1}^{\ast}$ since $\pi_{11}(L_1)\neq 0$. By
Lemma 2.2(b),  $L_1$ is an irreducible $L_0$-module of dimension
$2ml$. Therefore  dimension of $\pi_{11}(L_1)$ is equal to $2ml$
since $\pi_{11}(L_1)$ is irreducible $L_0$-module. On the other
hand, we have  $$(\dim{V_{1}})( \dim{W_{1}}) = \dim{ (V_{1}^{\ast}
\otimes W_{1}}) = \dim{\pi_{11}(L_1)} = 2ml.$$

Since  $V_{1}$ is a nontrivial $sp(m)$-module and $W_{1}$ is a
nontrivial $o(2l)$-module, it follows that  $\dim{V_{1}} \geq m$
and $\dim{W_{1}} \geq 2l$. Therefore  $\dim{V_{1}} = m$ and
$\dim{W_{1}} = 2l$. Hence $I_2$-module $W_{j_0}$ is standard.
\hfill $\fbox{}$

\end{proof}

$\quad$

\begin{lemma}

Let $S=K+L$ where ${ S} \cong {osp}(m,2n)$, ${ K} \cong osp(p,2n)$
and ${ L}\cong osp(m,2l)$. Then  for any $j \in \{1 \ldots d\}$,
$L_0$-module $W_{j}$ is not of the type 3.

\end{lemma}

\begin{proof}

We consider $L_0$-modules $W_1 \oplus \ldots \oplus W_d$. By
Lemmas 2.25 and  2.26, $L_0$-module $W_j$ is not of the type 1, 2
and 4. Hence any $L_0$-module $W_j$ is of the type 3. Moreover, by
Lemma 2.27, $L_0$-module $W_j$ has  dimension $2l$. Since since
$2l<2n$, it follows that $\dim{W_{j}} < \dim{W}$. Therefore $W$
contains at least two $L_0$-modules $W_1$ and $W_2$ of type 3. By
Corollary 2.18, $V$ is  an irreducible $I_1$-module. We consider
$I_1 \oplus I_2$-module $V^{\ast}_1 \otimes W_1$.

At first we prove that there exist $\lambda \in F$ such that
$\varrho_{12}(x)= \lambda\varrho_{11}(x)$ for any $x \in L_1$ (see
section 2.2). If $\pi_{12}(L_1)= \{0\}$ we choose $\lambda =0$.
Let $\pi_{12}(L_1)\neq 0$. Since $L_0$-module $L_1$ is
irreducible, $L_0$-module $L_1$ isomorphic to $L_0$-modules
$\pi_{11}(L_1)$ and $\pi_{12}(L_1)$. Hence   $\pi_{11}(L_1)$ and
$\pi_{12}(L_1)$ are isomorphic as $L_0$-modules. Therefore
$\varrho_{11}(L_1)$ and $\varrho_{12}(L_1)$ are also isomorphic as
$L_0$-modules. By Schur's Lemma,   $\varrho_{12}(x)=
\lambda\varrho_{11}(x)$ for any $x \in L_1$

Hence the matrix realization of $L_1$ has the  form
$$ \left\{
 \left(\begin{array}
                   { c  |  c  c  cc }
                     0  &  N_{11} & N_{12} & \ldots &  N_{1d}    \\ \hline
                     M_{11}   &   &        &        &    \\
                     M_{12}&   &        &        &    \\
                      \vdots  &   &        &   {\mbox{\Large 0}}     &   \\
                     M_{1d}&  &         &        &   \\

           \end{array}
   \right) \right\}
$$ where $M_{1i}$ is a  matrix  of order $2l \times m$,
$N_{1j}$ is a  matrix  of order $m \times 2l$ and $M_{12}= \lambda
M_{11}.$

The commutator of any two matrices from $L_1$ has the form:
$$ \left\{
 \left(\begin{array}
                   { c  |  c  c  cc }
                     X  &  0      & \ldots &  0         \\ \hline
                     0  &  Z_{11} & \ldots &  Z_{1d}    \\
                \vdots  &  \vdots &        &  \vdots    \\
                     0  &  Z_{d1} & \ldots &  Z_{dd}    \\
           \end{array}
   \right) \right\}
$$ where $Z_{ij}=M_{ji}N'_{jj}-M'_{ji}N_{jj}$.

We know that $L_0$-modules $W_1$ and $W_2$ are of type 3. Hence
there exist $Z_{11}\neq 0$ and $Z_{22}\neq 0$. Since
$Z_{11}=M_{11}N'_{11}-M'_{11}N_{11}$ and
$Z_{22}=M_{12}N'_{12}-M'_{12}N_{12}$, it follows that there exist
$M_{11}\neq 0$ and $M_{12}\neq 0$. Therefore $\lambda \neq 0$. On
the other hand, $Z_{21}=M_{12}N'_{11}-M'_{12}N_{11}=
\lambda(M_{11}N'_{11}-M'_{11}N_{11})= \lambda Z_{11}$.
 Hence $Z_{12} = \lambda Z_{11} \neq 0$. This contradicts the fact that $L_0$ has the form (16) where $Z_1 \neq 0$. \hfill $\fbox{}$

\end{proof}

\begin{corollary}
A Lie superalgebra  $ S \cong {osp}(m,2n)$ cannot  be decomposed
into the sum of two proper simple subalgebras ${ K}$ and ${ L}$ of
the  type $osp(p,2n)$ and $osp(m,2l)$, respectively.

\end{corollary}

\begin{proof}
The proof follows from the fact that for any $j \in \{1 \ldots
d\}$, $L_0$-module $W_{j}$ is not of the type 1, 2, 3 and 4.
\hfill $\fbox{}$
\end{proof}

$\quad$$\quad$

From Lemma  2.16 and Corollaries 2.24, 2.29, we obtain following
theorem
\begin{theorem}
A Lie superalgebra  $ S \cong {osp}(2k,2n)$ cannot  be decomposed
into the sum of two proper simple subalgebras ${ K}$ and ${ L}$ of
the  type $osp(p,2q)$ and ${osp}(s,2l)$, respectively.
\end{theorem}

From Theorems 2.14,2.15 and 2.30 we obtain  following Theorem:

\begin{theorem}

Let  $ S = {osp}(m,2n)$ be a  Lie superalgebra such that $S=K+L$
where $K$, $L$ are two proper basic simple  subalgebras.  Then $m$
is even, $m=2k$ and $K \cong osp(2k-1,2n)$, $L \cong sl(k,n)$ .
\end{theorem}

\end{document}